\documentclass[12pt,reqno]{article}
\usepackage{amsfonts}
\usepackage{graphicx}
\usepackage{amsmath}
\usepackage{amsthm}
\usepackage{latexsym,bm}
\usepackage{amssymb,float}
\usepackage{indentfirst}
\usepackage{mathrsfs}
\usepackage{caption}
 \usepackage{setspace}
\numberwithin{equation}{section}

\usepackage{hyperref}
\hypersetup{pdfnewwindow=true,colorlinks=true,linkcolor=blue,citecolor=blue,
	filecolor=magenta,urlcolor=blue 
	}

\begin{document}

\title{{\bf 
On Maur\'icio  M. Peixoto and the arrival of Structural 
Stability  to  Rio de Janeiro,  1955.
}}

\author{{\bf  Jorge Sotomayor\footnote{The author had the partial support of CNPq, PQ-SR - 307690/2016-4.
}
}}

\date{ }
\maketitle
\abstract{ This essay is an inquiry  about the circumstances, and subsequent   mathematical  consequences, of  the
  encounter, in 1955,  of Maur\'icio M. Peixoto (1921 - 2019) 
  with the work  of Henry F.  DeBaggis (1916 - 2002)  on   Structural Stability, Princeton 1952, developing  a notion introduced by
   Alexander A. Andronov and Lev S. Pontrjagin,  Gorkii  1937. 
          
  }

\section{ An accidental library encounter}\label{sec_library}

At the quiet  Library of the Institute of Pure and Applied Mathematics (IMPA), in  his  after  coffee  afternoon 
browsing books routine, 
Mauricio M. Peixoto (Peixoto), stumbled on  the article by H. F. DeBaggis (DeBaggis) entitled Dynamical Systems with Stable Structure, \cite{1952_b}, which captured his interest.  

He 
read attentively the  following phrases:

\vspace{0.3cm}

{\it 
$``\cdots $
For physical systems to perform certain operations they must, if they are to be useful, possess a certain degree of stability.
 Small perturbations should not affect the essential features of the system. 
 
 Since the physical components of a
system can never be duplicated exactly, experimental verification would be impossible unless the system remains stable under small variations. 
 
 The stability requirements of experiments provide a clue to the restrictions a mathematician should place on his nonlinear problems."  
 }

\vspace{0.3cm}
As an engineer, self taught mathematician and full professor of Rational Mechanics, with much higher 
 than average knowledge of the differential equations of Classical Mechanics, Peixoto  immediately noticed that a clever form  of stability -- or continuity-- was involved in 
 the author's text.  
  It was {\it Structural Stability}.
 
He glimpsed at the references and there he found the book  
{\it Theory of Oscillations}  of Andronov and Chaikin, 1949,  in the English Princeton translation  \cite{1949_ac},  of the Russian original \cite{1937_ac}, 1937. 
He had perused it in 1951 when 
writing his Thesis  
{\it The General Equations of Mechanics,} \cite{1951_b},  
presented in the competition for the  Chair of Rational Mechanics at the University of Brazil. 
The book  is cited in the  chapter {\it Stability}, of the Thesis, which deals  with the stability of stationary states, in the sense of Lyapunov, and not in that of the phase space depending,  as a whole,   on the analytic data defining the system.

On 1951, however, the notion of Structural Stability,  as mentioned along the Andronov and Chaikin book, did not attract his attention. The results   by Andronov and Pontrjagin, \cite{1937_ap},  stating a characterization of 
Structurally Stable Systems in a planar region, whose border they cross transversally, 
was 
formulated on an appendix of \cite{1949_ac}. No proofs were provided. 
This  omission  might have contributed to diminish the potential attractiveness of Structural Stability  for many mathematically inclined readers in the West.   
Two exceptions -- DeBaggis and Lefschetz -- will be discussed in Sections \ref{sec_baggis} and \ref{sec_lef}.

\vspace{0.3cm}

Complementing the  phrases above, in  a second reading round of \cite{1952_b},  Peixoto also glanced at:   
 
\vspace{0.3cm}
{\it`` $\cdots $
In this paper we give a complete treatment of the theory of structural stability. We have relaxed the conditions of analyticity which were imposed on the functions by \cite{1937_ap} and merely require that they have continuous first partials. " 
} 

\vspace{0.3cm}
 
 He focused again  on the concise list of references and attempted to  grasp the structure of the paper examining  the sections and  the  figures. 
  
In the last browsing round he also looked at the short final remarks about bifurcations --the lost of structural stability-- and the insinuation  of a possible continuation in that direction in  future work by the author.  

Time had elapsed  rapidly.  Approaching the librarian desk for borrowing the journal,   
he whispered:      
{\it This seems to be  good stuff.}

As he headed  back home,  carrying   in his briefcase the article of DeBaggis \cite{1952_b}, 
his thinking was accelerated. 
The  encounter, apparently 
 accidental,  followed by the 
active 
browsing 
through it, 
  had 
launched  Peixoto's  insightful  mathematical  intuition.

At home, he worked on  the article until late hours.  
Afterwards he also   discussed it with his wife Mar\'ilia C. Peixoto (Mar\'ilia), 1921 - 1961, also an engineer,  mathematician and his assistant in the Chair of Mechanics.

They  decided that, in due time,  he should  prepare a presentation on Structural Stability  
 
in the Seminar he directed, which had the participation of  Mar\'ilia, selected students and teaching  assistants of the Chair of Mechanics\footnote{Mentioned in Peixoto's  interview {\it Em\'eritos III, 2011},  
 \url{https://youtu.be/PToAegfcFKA.}}.

During the preparation of the Seminar, confusing points  in  \cite{1952_b}  were detected, where corrections  and  essential improvements were foreseen. The idea of the handy introduction of a Functional Space of Differential Equations, emerged. On it new results could be formulated  in topological terms.  This appears in the research  note  \cite{1955_abc} published in   September 1955, containing  innovative  results, such as the openness of structural stability.

Along the next months, the project of  a   visit to  Princeton to confer with  Solomon  Lefschetz, Editor of the Series Contributions
to the Theory of Nonlinear Oscillations,  began to be conceived.

\section {H. F. DeBaggis}\label{sec_baggis}
	In 1949,  on leave of absence form the University of Notre Dame  at Princeton, 
	Father Henry F. DeBaggis, CSC,  browsing through the recently published  book  by   Andronov and Chaikin \cite{1949_ac}, 
	encountered  the appendix on Structural Stability. 
	He was looking for new horizons to go beyond Axiomatic Hyperbolic  Geometry, field  in which he had written 
		a M.S. Dissertation and a PhD Thesis at   Notre Dame.
	
He decided to undertake the challenge  of giving a   proof of the  Theorem of Andronov and Pontrjagin stated there.  
Though, at the beginning,  he might  not have had the suitable  adroitness  in Differential Equations and Classical Analysis, Solomon Lefschetz supported him. DeBaggis participated in Lefschetz' 
Research Project  for two years.  See  Section \ref{sec_lef}.  

No  additional  information concerning DeBaggis'  research in mathematics was found. 
In his Obituary one can read that he taught Mathematics in several  North American Universities and also in other continents.   No chronological data  about him after he left Princeton  was found.

\section{Solomon Lefschetz, 1884 - 1972.  \label{sec_lef} }

We quote from Griffiths, Spencer and Whitehead  \cite{1992_nas}:

``{\it In 1943 Lefschetz became a consultant for the U.S. Navy
at the David Taylor Model Basin near Washington, D.C. There he met and worked with Nicholas Minorsky, 
who was a specialist on guidance systems and the stability of ships and who brought to Lefschetz' attention the importance
 of the applications of the geometric theory of ordinary differential equations to control theory and nonlinear mechanics. 
 
 From 1943 to the end of his life, Lefschetz' main interest was centered around ordinary nonlinear differential equations and their application
  to controls and the structural stabilities of systems.
  
   Lefschetz was almost sixty years old when he turned to differential equations, yet he did original work and stimulated research in this field as a gifted scientific administrator.''}
   
\vspace{0.03cm}
We quote from Lefschetz \cite{1959_onr},   p. 1: 

     {\it ``
    During most of World War II,
the undersigned, a consultant at the David Taylor Model Basin, had
frequent interviews with Dr. Nicolas Minorsky, in connection with
the latter's production of his well known 
Introduction to Nonlinear
Mechanics. 

Dr. Minorsky voiced repeated regrets at the impossibility
of creating in this country anything resembling the well known Institute
of Oscillations of Moscow, with its large staff of highly competent
mathematicians and physicists devoted to the problems of oscillations
and more generally to nonlinear mechanics (equally called nonlinear
differential equations). "
    }
    
    \vspace{0.03cm}
    
    We quote from Lefschetz  \cite{1959_onr}, p. 15: 
    
  {\it  `` Father Henry  DeBaggis, an 
Assistant 
 Professor of Mathematics
at Notre Dame and a Ph. D. under Karl Menger at Notre Dane also,
joined our Project in 1949 and was a member for two years 1949 - 51.

While his thesis was on 
Hyperbolic 
Geometry, he had little taste
for that subject and wished to change over to Differential Equations.\\
A reading of an appendix in Andronov-Chaikin \cite{1949_ac}
recently appeared, awoke his interest in structural stability. \\

This highly interesting concept had been launched in a Note of the Doklady
\cite{1937_ap} by Andronov and Pontrjagin. They considered a planar system
defined in a closed two-cell with 
vectors 
pointing outwards along
the boundary and asked under what conditions does the topology 
of the system of paths remain unchanged for small variations of
the 
vector field 
throughout the two-cell. \\
They stated n.a.s.c. for
this to happen but gave no proofs. 

\vspace{0.2cm}
DeBaggis undertook to establish
a complete theory and this objective was attained. 
His results were
developed in a paper which appear in
 \cite{1952_b}
and was sub
sequently translated into Russian \cite{1955_u}. 

It is only fair to say that DeBaggis derived very great benefit from discussions with M. Schiffer,  D. Spencer and the Director."\\
}

\vspace{0.03cm}  
   Lefschetz'  interest on Structural Stability was notorious. 
    Besides renaming it from  ``Roughness", he  had addressed  to Andronov a letter suggesting that  {\it ``a theorem of such importance should have its proof published."}   Unpublished correspondence  exhibited  in a glass display case,   together with manuscripts and books  of Andronov,   at the Meeting in commemoration of his  $100^{th} $  anniversary  \cite{2002_a100}.

 \section {The letter of Peixoto to Lefschetz, 1956.} 

Peixoto  proposed to establish that 
the systems defined by the conditions of Andronov and Pontrjagin is an open and dense set inside the open set of vector fields transversal to the border of the region supporting the system.

Upon reception of the letter, Lefschetz became deeply interested in Peixoto's project. 
He promptly wrote him back and suggested  how to proceed  to apply for a Research  Associate visiting position at Princeton.  

With the support of CNPq he visited Princeton during the academic year 1957 - 58, continuing further in RIAS, Baltimore, during 1958-59. 

This interesting stage of Peixoto mathematical activities, {\it The Golden Days}, in which he met Stephen Smale,   are recounted in his charming Speech of Acceptance of the TWAS Mathematics Award, 1987. See \cite{1987_peixoto_twas},  \cite{2019_peixoto_mu}.

\section{A meaning for  an   accidental encounter}

In 1943 the eminent mathematician Solomon Lefschetz:  algebraic geometer and topologist,  
reaching his sixties, 
decided  
 to  change his research interests and started to work in  Differential Equations. 
This happened after  encounters he had  with the  
Naval  
 Engineer  Nicholas Minorsky, an  {\it ``ingeni\'eur savant"}  \cite{2015_petit}. See Section \ref{sec_lef}.
 
From Princeton Lefschetz launched  several editorial activities to  spread the interest on Nonlinear Differential Equations, \cite{1959_onr}. 
One of them was the Edition of the Series  {\it Contributions to Nonlinear Oscillations}, which in Vol. II   published the  paper of   DeBaggis in 1952.  In 1955   it    reached Peixoto in Rio de Janeiro and,  
after  substantial  elaboration,  took him to Princeton to  have the benefit of  {\it The Golden Days} of his mathematical career.   

{\verse {\it

On the afternoon of the encounter of Peixoto with\\
 DeBaggis paper, a narrow cosmic 
 fracture
  allowed \\
   a considerable amount of repressed 
mathematical energy to\\ flow  into the library of IMPA.\\
It was liberated by the  apparently accidental and seemingly\\ naive
 ceremony of browsing a mathematics article,\\
as a dam under pressure releases  its flow through a fissure. \\
On  the right side of the dam 
to channel the flow, \\
 there was a receptive,  insightful and sensitive  mind. \\
 One 
 that was   prepared for the challenge  and could grasp the \\complexity of the  
 inflowing 
 mathematical
  message.\\
  } }

  \vspace{0.5cm}

\noindent 
  Thus, Structural Stability, carrying  its  rich heritage and 
  mathematical background, arrived from Gorkii, now  Nizhnii Novgorod,   to Rio de Janeiro, then the Capital of Brazil, 
   after an  auspicious  maturation   in Princeton, one of the capitals of World Mathematics.      
   
   In 1987 Peixoto received the TWAS Mathematics Award for  his contribution to 
   Structural Stability in Dynamical Systems  \cite{1987_peixoto_twas}. 
   
 \section{Contacts and influences of Peixoto around Princeton 1958,  TWAS 1987 and after }
 
\begin{figure}[H]
\begin{center}
\includegraphics[scale=1.00]{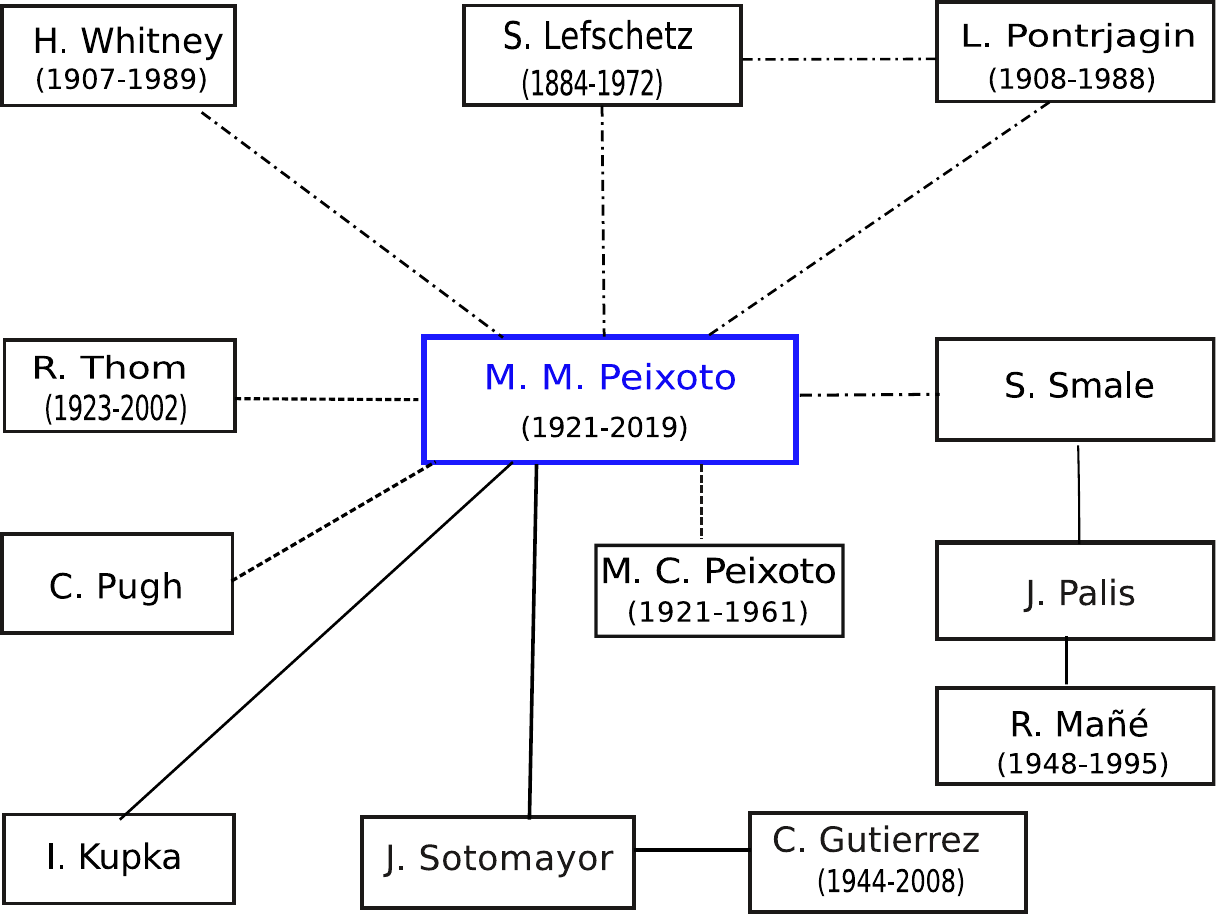}

 \begin{minipage}[t][3cm][b]{0,7\textwidth}  
{\small  Mathematical contacts and influences mentioned by Peixoto  in his charming   Acceptance Speech on the occasion  of the TWAS Mathematics Award that he received in 1987. There he presented an outline of crucial developments in  Dynamical Systems.  See \cite{1987_peixoto_twas} which had a Portuguese translation included in the obituary essay \cite{2019_peixoto_mu}. }
	\end{minipage}    
  
\end{center}
\end{figure}

Other aspects of Peixoto's mathematical biography can be found in \cite{2019_sotomayor_list}, 
\cite{2018_sotomayor_encounters} and \cite{2001_sotomayor_peixoto80}.

  \vspace{0.5cm}

 \noindent {\bf Acknowledgement.}  The author is grateful to  J. Derwent  and to H. Weiss  for their valuable  help to get access to, respectively,  the works of DeBaggis  and   the enlightening Lefschetz Report \cite{1959_onr}. 
 Thanks are also due to D. Schlomiuk,  M. Sotomayor,  L. F. Mello,  A. R. da Silva and R. A. Garcia for their  helpful comments to this text.\\

\vspace{0.4cm}
\noindent Instituto de Matem\'atica e Estat\'istica, \\Univ. de S\~ao Paulo, Rua do Mat\~ao, 1010,\\ CEP: 05508 - 090, \\S\~ao  Paulo, S.P, Brazil\\
E-mail address: sotp@ime.usp.br

\end{document}